\documentclass{amsart}[12pt]
\usepackage{amsmath,amssymb,amsthm,amsfonts,graphicx}
\usepackage[dvipsnames]{xcolor}
\usepackage[notref,notcite]{showkeys}  

\usepackage{amssymb}

\newtheorem{them}{\bf{Theorem}}

\newtheorem{prop}{Proposition}

\renewcommand{\theprop}{\arabic{prop}}
\newtheorem{lem}[prop]{Lemma}

\newenvironment{pf}{{\noindent \bf Proof. \quad}}{\hfill $\Box$}

\renewcommand{\theequation}{\arabic{equation}}

\newcommand{\be}{\begin{equation}}
\newcommand{\ee}{\end{equation}}
\newcommand\bes{\begin{eqnarray}} \newcommand\ees{\end{eqnarray}}
\newcommand{\bess}{\begin{eqnarray*}}
\newcommand{\eess}{\end{eqnarray*}}

\def\Xint#1{\mathchoice
  {\XXint\displaystyle\textstyle{#1}}%
  {\XXint\textstyle\scriptstyle{#1}}%
  {\XXint\scriptstyle\scriptscriptstyle{#1}}%
  {\XXint\scriptscriptstyle\scriptscriptstyle{#1}}%
  \!\int}
\def\XXint#1#2#3{{\setbox0=\hbox{$#1{#2#3}{\int}$}
  \vcenter{\hbox{$#2#3$}}\kern-.5\wd0}}

\def\dashint{\Xint-}

\begin{document}



\title[Evolution of Locally Convex Closed Curves in Nonlocal Curvature Flows]{Evolution of Locally Convex Closed Curves in the Area-Preserving and Length-Preserving Curvature Flows}

 \author[N. Sesum]{Natasa Sesum}
 \author[D.H. Tsai]{Dong-Ho Tsai}
 \author[X.L. Wang]{Xiao-Liu Wang}
 \address[N. Sesum]{Department of Mathematics, Rutgers University, Pitscataway 08854, USA.\ \ \ \ \ \ \ \ \ \ \ \  E-mail: \ \textit{natasas@math.rutgers.edu}}
 \address[D.H. Tsai]{Department of Mathematics, National Tsing Hua University, Hsinchu 300, TAIWAN. E-mail:\ \textit{dhtsai@math.nthu.edu.tw}}
 \address[X.L. Wang]{School of Mathematics, Southeast University, Nanjing 211189, PR CHINA. E-mail: \ \textit{xlwang@seu.edu.cn}}


\maketitle

\begin{abstract}
 We provide sufficient conditions on an initial curve for the area preserving and the length preserving curvature flows of curves in a plane, to develop a singularity at some finite time or converge to an $m$-fold circle as time goes to infinity.  For the area-preserving flow, the positivity of the enclosed algebraic area determines whether the curvature blows up in finite time or not, while for the length-preserving flow, it is the positivity of an energy associated with initial curve that plays such a role.
\end{abstract}




\renewcommand{\theequation}{\thesection.\arabic{equation}}
\setcounter{equation}{0} \setcounter{prop}{0}
\renewcommand{\theprop}{\arabic{section}.\arabic{prop}}
\section{Introduction}

\subsection{Background}

The planar curvature flows, arising in many application fields,
such as phase  transitions, image processing, etc., have received
a lot of attention in  recent years. Generally, their evolution equations
take the form of
 \begin{equation}\label{eqn1}
 \left\{
\begin{array}{ll}
{\partial X(u,t)}/{\partial t}= F {\mathbf{n}},\\
 X(u,0)=X_{0}(u),
\end{array}\right.
\end{equation}
where $X(u,t):S_m^1\times[0,T)\rightarrow\ {\mathbb{R}}^2\ (T>0)$ is
a family of evolving curves with speed $F$ along inward pointing
unit normal ${\mathbf{n}}$ and $X_0$ is a
 closed curve with total curvature of $2m\pi\ (m\in
\mathbb{Z}^+)$. When $F$ is the signed curvature $\kappa(u,t)$ at $X(u,t)$,
(\ref{eqn1}) is the famous curve shortening flow, which evolves an embedded closed curve into a convex one and then shrinks it
into a round point (see for example the pioneering works of Gage-Hamilton \cite{GH1} and Grayson \cite{Gr1}). If $X_0$ is immersed
and locally convex \footnote{ Here and after, we use the convention that for locally convex
plane curves the curvature  is positive everywhere.}, the behaviour of the curve shortening flow becomes more complicated and
has been studied by Abresch-Langer \cite{AL1} and Angenent \cite{Ang1}.
When $F = \frac{1}{\alpha}|\kappa|^{\alpha-1}\kappa$ (for $\alpha\neq 0$) and
$X_0$ is a locally
 convex curve (including the embedded or immersed case), (\ref{eqn1}) is known as the generalized or power-type curvature flow.
The different homogeneous degree $\alpha$ of $F$
 w.r.t. the curvature $\kappa$ could possibly result into different evolution behavior of the flow,
 see for instance \cite{And1, And2,Ur1} and etc.
 For more about (\ref{eqn1}) and its applications, one may refer to monographs \cite{Cao1,CZ1,Vi1} and references therein.

 Another class of popular curvature flows are the so-called  nonlocal curvature flows, whose
  evolution equations take the form of
 \begin{equation}\label{eqn2}
 \left\{
\begin{array}{ll}
{\partial X(u,t)}/{\partial t}=[f(\kappa(u,t))-\lambda(t)]{\mathbf{n}},\\
 X(u,0)=X_{0}(u).
\end{array}\right.
\end{equation}
 As for the speed, $f(\kappa)$ is a given function of
curvature satisfying $f'(z)>0$ for all $z$ in its domain, and
$\lambda(t)$ is a function of time which may depend on certain
global quantities of $X(.,t)$, say enclosed algebraic area $A(t)$, length
$L(t)$, or others. Specifically, $A(t) = \int_{R^2}w(x,y,t)dxdy$ and
$w(x,y,t)$ is the winding number of $\gamma(\theta,t)$ around
$(x,y)\in R^2$. When these flows are applied to the image processing, they have better boundary smoothing effect than the curve shortening flow, see
\cite{ST1,XY1}.  Also, some nonlocal flows could be used to describe the motion of the
interface arising in nonlocal models of Allen-Cahn equation or Hele-Shaw models, see \cite{CFM1} and \cite{DM1} respectively.

\vskip 10 pt

The purpose of this paper is to study $\kappa^{\alpha}$-type nonlocal flow
(\ref{eqn2}) with the speed function given by

\bess
\ \ \ \ \ \ \ \ \ \ \ \ \ \ \ \ f\left(  k\right)  -\lambda\left(  t\right)  =\kappa^{\alpha}-\frac{1}{{L}\left(
t\right)  }\int_{X\left(  \cdot,t\right)  }\kappa^{\alpha}ds{,\ \ \ \alpha
> 0,\ \ \ } \ \ \ \ \ \ \ \ \ \ \ \ \ (AP)\nonumber%
\eess
or$\ $%
\bess
\ \ \ \ \ \ \ \ \ \ \ \ \ \ f\left(  k\right)  -{\lambda}\left(  t\right)  =\kappa^{\alpha}-\frac{1}{{2m\pi}}%
\int_{X\left(  \cdot,t\right)  }k^{\alpha+1}ds{,\ \ \ \alpha> 0,\ \ \ }%
\ \ \ \ \ \ \ \ \ \ \ \ \ (LP)\nonumber%
\eess
and  initial curve $X_0$ being smooth, immersed, locally convex and closed.
Here $2m\pi(m\in\mathbb{Z}^+)$ denotes the total curvature of $X(\cdot,t)$, $s$ is the arc length parameter
of $X(\cdot,t)$ and the constant $\alpha>0$ is arbitrary.  We shall see shortly
that under the evolution equation (AP) the flow is area-preserving, and under the equation (LP) the flow is length-preserving.  The abbreviations
AP and LP will be used to indicate that the flow is ``area-preserving" and ``length-preserving", respectively. Without causing ambiguity in the context, we define
$$
\displaystyle{\lambda(t) = \frac{\int_X
\kappa^\alpha\,ds}{\int_X
\,ds} = \frac{\int_X
\kappa^\alpha\,ds}{L(t)}}$$
for the AP flow and
$$
\displaystyle{\lambda(t) = \frac{\int_X
k^{\alpha+1}\,ds}{\int_X
k\,ds} = \frac{\int_X
k^{\alpha+1}\,ds}{2m\pi}}$$
for the LP flow.

When initial curve
$X_0$ is embedded, convex and closed, Gage studied the AP flow with $\alpha=1$  in \cite{Ga2} and  showed  the flow could exist for all time,
preserving the evolving curves' enclosed area, while
decreasing their length, and finally making them converge to a round circle in
$C^\infty$ metric. After that, many authors had nice contributions to the research of planar nonlocal curvature flows. For example, Ma-Zhu \cite{MZ1} studied the
LP flow with $\alpha=1$, and  Jiang-Pan \cite{JP1} studied the gradient flow of isoperimetric ratio functional, which increases the enclosed area of evolving curves and decreases their length. Those flows also exist globally and converge smoothly to round circles. Recently, the authors considered the AP flow and the LP flow with any $\alpha>0$ and other nonlocal flows in \cite{TW1}, showing that all of those flows have the same convergence behavior as in \cite{Ga2}. See also related work \cite{MPW1} by Pan et al. The key step in proving the above flows' convergence is showing the time-independent bounds for the curvature of evolving curves, which is obtained in \cite{Ga2,MZ1,JP1,MPW1}, by modifying the proof of Gage-Hamilton \cite{GH1}, and in \cite{TW1} by using the support function method of Tso \cite{Tso1}. In any case, the Bonnesen inequality plays an important role.

When initial curve $X_0$ is immersed, locally convex and closed, it is not hard to show
that the evolving curves' enclosed algebraic area (length) is preserved under the AP
flow  (the LP flow respectively), but some things become more difficult and different from a case of embedded convex and closed curves.
For instance, it is unknown whether Bonnesen inequality holds or not in this case.
This requires  developing new methods when dealing with the evolution of these flows.  In
\cite{EI1}, among other things, Escher-Ito showed that the singularity must
happen (that is, the
curvature  blows up) at some finite time in the AP flow with $\alpha=1$ when the algebraic area of initial curve is negative.
Later, in \cite{WK1}, two
classes of rotationally symmetric and locally convex initial curves, namely, highly symmetric curves and
 Abresch-Langer type curves (see the definitions in Section 1.2), both enclosing positive
algebraic area, are
found to guarantee the convergence of the AP flow with $\alpha=1$
 to $m$-fold circles. A similar result is established for the LP flow with $\alpha=1$ in \cite{WW1}.

 \vskip 10 pt

In this paper, we would like to investigate the AP flow and the LP flow for any $\alpha>0$ for immersed, locally convex, closed curves. The previous results about singularity formation and convergence when $\alpha=1$ (as discussed in \cite{EI1,WK1,WW1}) are generalized to the case when $\alpha > 0$. Moreover, by observing the sufficient conditions on finite-time singularity or global convergence of the flow, we can compare the difference between the AP flow and the LP flow, and also the difference between the nonlocal flows and the curve shortening flows. The key ingredient is resolving the convergence problem for a globally existing flow. To do this, one faces a traditional problem in the theory of nonlinear evolution equations to prove that under some hypotheses, any global solution to a particular nonlinear evolution equation is uniformly bounded from above. The arguments are motivated by the works of Chou \cite{Ch1} and  Dziuk-Kuwert-Sch$\ddot{\mathrm{a}}$tzle \cite{DKS1}.


\vskip 10 pt

Before ending this section, we would like to mention other nonlocal flows. In the higher dimensional Euclidean space, people consider nonlocal flows for hypersurfaces.  For example, there are
volume preserving and surface area preserving mean curvature flows for embedded closed convex hypersurfaces, see Huisken \cite{Hui1}, McCoy \cite{Mc1} and etc. When the speed function is a nonlinear function of curvature, see a recent work \cite{Si1}.  People also considered the problems of nonlocal flow with boundary, see \cite{Mad1,SK1}, for the planar curve case and \cite{AK1,CM2} for the higher dimensional cases.




\vskip 10 pt

\subsection{Main Theorems}

It turns out that both flows preserve local convexity during the evolution, which we prove in  Lemma \ref{convexity}. Therefore, we can always use the normal angle $\theta\in[0,2m\pi]$ for the parametrization of the evolving curves. Here and after, we use $I$
 to denote the circle:
$$I = [0,2m\pi],$$
or more precisely
$$I ={\mathbb{R}}/2m\pi{\mathbb{Z}}.$$

We first introduce two classes of rotationally symmetric curves.
  The first class is highly symmetric curves, which are defined to be
  $${\mathcal{H}}_{m,n} = \{The \ locally \ convex \ curves \ with \
   n-fold \ rotational$$
   $$ \ \ \ \ \ \ \ symmetry \ and \ total \ curvature \ of\  2m\pi \ (m\ and \ n  $$
   $$ are \ mutually \ prime \ and \ satisfy \ n>2m)\}.$$
The second class is Abresch-Langer type curves.
Before we explain what those curves are, let us define the support function of a locally convex closed curve $X$
(parameterized by its normal angle $\theta\in I$) to be,
\begin{equation}
\label{eq-support-fun}
h(\theta) = <X(\theta),-{\mathbf{n}}(\theta)>,
\end{equation}
where $-{\mathbf{n}}$ is the outward normal vector. The relationship
between locally convex, closed curves and their support functions is
contained in Proposition 2.1 of \cite{CZ1}.
Now we define the Abresch-Langer type curves to be
$${\mathcal{A}}_{m,n} = \{
The \ locally \ convex \ smooth \ curves \ having \ n-fold \
rotational $$
$$\ \ \ \ \ \ \ \ \ symmetry \ and \ total \ curvature \ of \ 2m\pi \ (m \
and \ n \ are \ mutually$$
$$\ \ \ \ \ \ \ \ \ prime \ and \ satisfy \ n<2m),\ and \ having\ the \ property \ (\mathcal{P})\},$$
where the property ($\mathcal{P}$) means

\vskip5 pt

($\mathcal{P}$) Under the normal angle's parametrization,  the locally convex curve's support function $h(\theta)$ and curvature function $\kappa(\theta)$ are symmetric with respect to $\theta = 0$ and $\theta = m\pi/n$; both of them are strictly
decreasing in $(0,m\pi/n)$; moreover, $h(m\pi/n)>0$.

\vskip 5 pt

\begin{figure}
\includegraphics[width=1.25\textwidth]{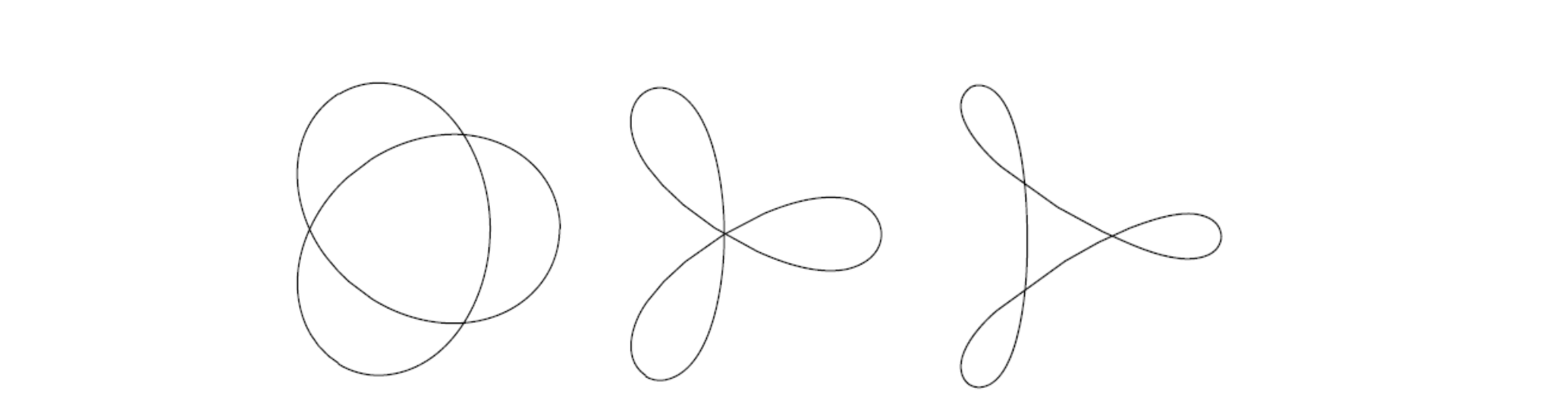}
\caption{The deformation of Abresh-Langer type curve }
\end{figure}

\begin{them}\label{them1}
Let the initial curve $X_0$ be immersed, locally convex and closed. Then the following holds for the AP flow starting at $X_0$.
\begin{enumerate}
\item
\label{part-one}
If $A_0 < 0$ or $L_0^2 < 4m\pi A_0$ (where $A_0$ and $L_0$ are the enclosed algebraic area and the length of $X_0$, respectively), then a  singularity occurs during the evolution of the AP flow.

\item
\label{part-two}
If $X_0$ encloses a zero algebraic area, that is,  $A_0 = 0$,  then a singularity appears at the maximal existence time $T_{\mathrm{max}}$. If $T_{\mathrm{max}}=\infty$, then the flow converges to a point.

\item
\label{part-three}
If $X_0\in \mathcal{H}_{m,n}$ with $n>2m$,  then the AP flow exists globally and
  converges to an $m$-fold circle in $C^\infty$-metric as time goes to infinity.
\item
\label{part-four}
If $X_0\in \mathcal{A}_{m,n}$ with $n<2m$,  then the AP flow exists globally and
  converges to an $m$-fold circle in $C^\infty$-metric as time goes to infinity.
 \end{enumerate}
\end{them}

\noindent \textbf{Remark 1.}  In \cite{EI1}, the authors proposed a question whether the maximal existence time of the AP flow with $A_0=0$  is finite or not? Our result in Theorem \ref{them1} (2) implies that if one can show a locally convex closed curve with $A_0=0$ does not evolve into a point as $t\rightarrow T_{\max}$, then $T_{\max}$ must be finite.

\vskip 5 pt

\noindent \textbf{Remark 2.}  An example of a curve in $\mathcal{H}_{m,n}$ is a pentagram. In the work \cite{EG1} of Epstein and Gage, all the curves   $\mathcal{H}_{m,n}$ are shown to have positive support function if the symmetric center is chosen to be the origin. The examples of curves that belong to $\mathcal{A}_{m,n}$ can be found in \cite{AL1, And2, Au1} (Abresch-Langer curves, the self-shrinkers in the curve shortening flow).
In Figure 1, one can see that an Abresch-Langer type curve could be deformed artificially into a curve with negative algebraic area.

\vskip 5 pt

\noindent \textbf{Remark 3.}  When $X_0$ is a rotationally symmetric locally convex curve with $m=1$ and $n\geq2$, it is just an embedded convex closed curve with  symmetry. In fact, the AP flow for general embedded convex closed curves has been investigated in \cite{Ga2} and \cite{TW1}.

\vskip 10 pt

 Our second main theorem about the LP flow is stated as follows. We define the following energy for the curvature function $\kappa$ parametrized by normal angle
$\theta\in I$,

$$E(t) = \int_I((\kappa^\alpha)_{\theta})^2\,d\theta - \int_I(\kappa^\alpha-\overline{\kappa^\alpha})^2\,d\theta,$$
where
$$\overline{\kappa^\alpha} = \frac{1}{2m\pi}\int_I\kappa^\alpha\,d\theta.$$

\begin{them}\label{them2}
Let the initial curve $X_0$ be immersed, locally convex and closed. Then the following holds for the LP flow starting at $X_0$.
\begin{enumerate}
\item If  $X_0$ satisfies that
$E(0)\leq 0$ and $k_0$ is nonconstant on $I$, or if $L_0^2 < 4m\pi A_0$ (where $A_0$ and $L_0$ are the enclosed algebraic area and the length of $X_0$, respectively), then a singularity appears during the evolution of
the LP flow;

\item If $X_0$ has $2m\pi$
  total curvature and $n$-fold rotational symmetry with $m/n\leq 1$, then the LP flow  exists globally and
  converges to an $m$-fold circle in $C^\infty$-metric as time goes to infinity.
  \end{enumerate}
\end{them}

\noindent \textbf{Remark 3.} We note that by the Poincar$\acute{\mathrm{e}}$ inequality (or Wirtinger inequality) any initial curve with period $|I|$ satisfies
$$\int_I(\kappa_0^\alpha-\overline{\kappa_0^\alpha})^2\,d\theta \leq \Big(\frac{|I|}{2\pi}\Big)^2\int_I((\kappa_0^\alpha)_{\theta})^2\,d\theta.$$
From this inequality, we can see that the curvature function of a rotationally symmetric curve $X_0$ in Theorem \ref{them2} (2)  satisfies $E(0)>0$ if $m/n<1$ and $E(0)\geq 0$ if $m=n=1$. The convergence result in Theorem \ref{them2} (2) coincides with the existing results in \cite{MZ1} and \cite{TW1} for the case of embedded convex curves, i.e., the case $m=1$.

\vskip 20 pt

We can see from the above theorems that there exists an essential difference between two nonlocal flows.
 For the AP flow, the positivity of the enclosed algebraic area determines whether the curvature  blows-up or not, while for the LP flow, it is the positivity of an energy associated with initial curve that plays such a role. Also, nonlocal flows exhibit different behavior from the (normalized) curve shortening flow. The homogeneous degree $\alpha$ of $F$  w.r.t. the curvature $\kappa$ does not affect the evolution behavior of considered nonlocal flows.

\vskip 10 pt

 This paper is organized in the following way. We reformulate our problems in Section 2 and give some basic lemmas in Section 3. Then we prove Theorems \ref{them1} and \ref{them2} separately
in Sections 4 and 5, respectively.

\renewcommand{\theequation}{\thesection.\arabic{equation}}
\setcounter{equation}{0} \setcounter{prop}{0}
\renewcommand{\theprop}{\arabic{section}.\arabic{prop}}
\section{Reformulation of problems}

For the convenience of readers, we use the following notations:

\vskip 10 pt

\ \ $ds$\ \ \ \ \ \ \ \,\, the differential element of arc-length

\ \ $\theta$\ \ \ \ \ \ \ \ \ \,  the normal angle of
$X(\cdot,t)$

\ \ ${\mathbf{n}}$\ \ \ \ \ \ \ \ \ \,  the inward pointing unit normal of
$X(\cdot,t)$

\ \ $L(t)$\ \ \ \ \ \ \ the length of $X(\cdot,t)$

\ \ $A(t)$\ \ \ \ \ \ \ the algebraic area of $X(\cdot,t)$ defined by
$-\frac{1}{2}\int_X<X,{\mathbf{n}}>\,ds$

\ \ $\kappa(\cdot, t)$\ \ \ \ \ \,the curvature of  $X(\cdot,t)$ w.r.t
${\mathbf{n}}$.

\ \ $h(\cdot,t)$ \ \ \ \,   the support function of curve $X(\cdot,t)$ as given in \eqref{eq-support-fun}.

\vskip 10 pt

\noindent Here, we always take the orientation of $X(\cdot,t)$ to be counter clockwise.

\vskip 10 pt

The evolution of various geometric quantities
along both flows  can be deduced from the general
formulas in \cite{CZ1} and \cite{Ga2}.

%
%

%

\vskip 5 pt

 \ \ \ $\displaystyle{\frac{\partial \kappa}{\partial t} = (\kappa^\alpha)_{ss} + \kappa^2(\kappa^\alpha -
 \lambda(t))},$

\vskip 5 pt

 \ \ \ $\displaystyle{\frac{d L}{d t} = - \int_X \kappa(\kappa^\alpha -
 \lambda(t))ds},$

\vskip 5 pt

 \ \ \ $\displaystyle{\frac{d A}{d t} = - \int_X (\kappa^\alpha -
 \lambda(t))ds}.$

\vskip 10pt

 \noindent Here, it can be easily seen that the enclosed algebraic
 area $A(t)$ of the evolving curves is preserved while the length $L(t)$ is non-increasing along the AP flow. And for the LP flow, $L(t)$ is preserved while $A(t)$ is non-decreasing. Indeed, we know (see \cite{Ga2}) that for a family of
time-dependent closed curves $X(\theta,t):{I}\times [0,T)\to
{\mathbb{R}}^2$ with time variation $\partial X(\theta,t)/\partial t
= W(\theta,t)\in {\mathbb{R}}^2$, their length $L(t)$ and enclosed
algebraic area $A(t)$ satisfy the following:
$$\frac{dL}{dt}(t) = -\int_X <W,\kappa {\mathbf{n}}>\,ds, \ \
\frac{dA}{dt}(t) = -\int_X <W,{\mathbf{n}}>\,ds.$$
For the AP flow, it is just
straightforward to verify $dA(t)/dt\equiv 0$ on $(0,T)$ and
$$\frac{dL}{dt}(t) = -\int_X \kappa\Big(\kappa^\alpha - \frac{\int_X\kappa^\alpha\,ds}{L(t)}\Big)\,ds.$$
To show
$dL(t)/dt\leq 0$, it is equivalent to show
$$\int_X\,ds\int_X \kappa^{\alpha+1}\,ds - \int_X \kappa\,ds\int_X\kappa^\alpha\,ds \geq 0.$$
Indeed, this holds by multiplying out the following H$\ddot{\mathrm{o}}$lder inequalities:
$$\int_X \kappa\,ds \leq \Big(\int_X \kappa^{\alpha+1}\,ds\Big)^{\frac{1}{\alpha+1}}\Big(\int_X \,ds\Big)^{\frac{\alpha}{\alpha+1}}$$
and
$$\int_X \kappa^\alpha\,ds \leq \Big(\int_X \kappa^{\alpha+1}\,ds\Big)^{\frac{\alpha}{\alpha+1}}\Big(\int_X \,ds\Big)^{\frac{1}{\alpha+1}}.$$
For the LP flow, the conclusion could be proved via a similar argument.

 \vskip 10 pt

When locally convex solution $X(\cdot,t)$ is considered,
 each point on it has a unique tangent and one can use the
normal angle $\theta \in I $ to
parameterize it. Generally speaking, $\theta$ is a function
depending on $t$. In order to make $\theta$ independent of time $t$,
one can attain that by adding a tangential component to the velocity
vector $\partial X/\partial t$, which does not affect the geometric
shape of the evolving curve (see, for instance, \cite{Ga2}). Then
the evolution equations can be expressed in the coordinates of $\theta$
and $t$. If we denote by
$\kappa(\theta,t)$ the curvature function of $X(\theta,t)$, Problem
 (AP) or (LP) could be reformulated equivalently
 as follows \bes
\left\{\begin{array}{ll} \kappa_{t}=\kappa^{2}[(\kappa^\alpha)_{\theta\theta} +
\kappa^\alpha - \lambda(t)],\ \
(\theta,t)\in {I}\times(0,T_{\max}),\\
\kappa(\theta,0) = k_0(\theta), \ \ \theta\in {I}
\end{array}
\right. \label{eqn2.1}\ees where $k_0$ is the curvature of
 $X_0$. In terms of the new coordinates, we have
 $$\lambda(t) = \frac{\int_I \kappa^{\alpha-1}(\theta,t)\,d\theta}{L(t)} $$
 for the AP flow,
 and
$$\lambda(t) = \frac{\int_I \kappa^\alpha(\theta,t)\,d\theta}{2m\pi}$$
for the LP flow.

\vskip 10 pt

The evolution equation of a support function $h(\theta,t)$, along the AP flow or the LP flow, is given by

 \bes
\left\{\begin{array}{ll} h_t =
-(h+h_{\theta\theta})^{-\alpha}+
\lambda(t), \ \ (\theta,t)\in {I}\times(0,T_{\max}),\\
h(\theta,0) = h_0(\theta), \ \ \theta\in {I},
\end{array}
\right. \label{eqn2.2}\ees
where $h_0$ is the support function of $X_0$.

\vskip 10 pt

For brevity, we set $v=\kappa^\alpha$ and $p=1+\frac{1}{\alpha}$. Then $v$ satisfies
\begin{equation}\label{curv_eqn_2}
v_t = \alpha v^p(v_{\theta\theta}+v-\lambda(t)), \  \  (\theta,t)\in I \times (0,T_{\max}).
\end{equation}
This equation will be frequently used later.

\renewcommand{\theequation}{\thesection.\arabic{equation}}
\setcounter{equation}{0} \setcounter{prop}{0}
\renewcommand{\theprop}{\arabic{section}.\arabic{prop}}
\section{Some basic lemmas}

In this section, we  prove a few lemmas that are needed to complete proofs of Theorems \ref{them1} and \ref{them2}. The first one is about the unique existence of both flows.

\begin{lem} When the initial curve is immersed, locally convex, closed and smooth, both flows have unique smooth, locally convex solutions on a time interval $[0,T)$. If $T_{\max}$  is the maximal existence time, then either $T_{\max} = \infty$ or $\limsup_{t\to T_{\max}} \max \kappa(\theta,t) = \infty$.
\end{lem}
\begin{pf}
The unique existence of the flow can be proven by applying the classical Leray-Schauder fixed point theory to Problem (\ref{eqn2.1}). See details in \cite{MPW1}, where a generalized area-preserving flow is studied. One can also find the relative references in \cite{Mc1}, where the nonlocal flows in higher dimensions are discussed. The preserved convexity will be proven in the next lemma.
\end{pf}


\vskip 5 pt

By the maximum principle, we can show that the local convexity of initial curve
 is preserved along both flows.
\begin{lem}\label{convexity}If the initial curve $X_0$ is locally convex, then $X(\cdot,t)$ continues to be locally convex as long as
the flow  exists.
\end{lem}
\begin{pf} By the continuity, $\min_{\theta\in I}\kappa(\theta,t)$ keeps positive on small
time interval. Assume that the time span of the flow is $T_{\max}$. Suppose that the conclusion is
not true. Then there must exist the first time, say $t_1<T_{\max}$, such that $\min_{\theta\in I}\kappa(\theta,t_1)=0$.
Next, consider the quantity
\bess
\Phi(\theta,t)=\frac{1}{\kappa(\theta,t)}
-\frac{L(t)}{2m\pi}-\frac{1}{2m\pi}\int_{0}^{t}\int_{0}^{2m\pi
}\kappa^\alpha(\theta,\tau)  d\theta d\tau
\eess
with $(\theta,t)\in I\times[0,t_1)$. By (\ref{eqn2.1}), we have
\bess
\Phi_{t}  & =& -\left(  \kappa^{\alpha}\right)  _{\theta\theta}-\kappa^{\alpha}%
+\lambda\left(  t\right)  -\left(  \lambda\left(  t\right)  -\frac{1}{2m\pi
}\int_{0}^{2m\pi}\kappa^{\alpha}d\theta\right) \nonumber \\
& & -\frac{1}{2m\pi}\int_{0}^{2m\pi
}\kappa^{\alpha}\left(  \theta,t\right)  d\theta\nonumber\\
&=& \alpha \kappa^{\alpha+1}\Phi_{\theta\theta}  -\alpha(\alpha+1)\kappa^{\alpha+2}
\Phi_{\theta}^{2}  - \kappa^{\alpha}\\
&\leq& \alpha \kappa^{\alpha+1}\Phi_{\theta\theta}.
\eess
Hence by the maximum principle,
\bess
\frac{1}{\kappa(\theta,t)}\leq\max_{\theta\in I}\left(\frac
{1}{\kappa_{0}(\theta)}\right)  +\frac{L(t)  -L(
0)}{2m\pi}+\frac{1}{2m\pi}\int_{0}^{t}\int_{0}^{2m\pi}\kappa^\alpha(
\theta,\tau)  d\theta d\tau
\eess
for all $(\theta,t)\in I\times[0,t_1)$. In fact, except the fact that $L(t)$ is nonincreasing in
both flows, we have the following observations:
$$
\max_{\theta\in I}\left(\frac{1}{\kappa_{0}(\theta)}\right)
-\frac{L(0)}{2m\pi}=\max_{\theta\in I}\left(\frac{1}
{\kappa_{0}(\theta)}\right) - \frac{1}{2m\pi}\int_{0}^{2m\pi}\frac
{1}{\kappa_{0}(\theta)}d\theta \geq 0,
$$
and
$$\max_{(\theta,t)\in I\times[0,t_1)} \kappa(\theta,t) \leq C_1(t_1) < \infty$$
 for some constant $C_1(t_1)$.
Therefore,
$$\min_{\theta\in I}\kappa(\theta,t) \geq C_2(t_1) > 0, \ \forall\ t\in[0,t_1),$$
for some constant $C_2(t_1)$. This contradicts our assumption that $\min_{\theta\in I}\kappa(\theta,t_1)=0$ and hence finishes the proof of the Lemma.
\end{pf}

\vskip 5 pt

In the following, we establish the gradient estimate for the curvature of evolving curves.
\begin{lem}
\label{grad_est}Under both flows, there holds the estimate%
\begin{equation}
\max_{I\times\left[  0,t\right]  }\Psi\leq\max\left\{  \max_{I%
\times\left[  0,t\right]  }v^{2},\ \ \max_{I\times\left\{  0\right\}
}\Psi\right\}  ,\ \ \ \forall\ t\in\lbrack0,T_{\max}), \label{v-grad-Ben}%
\end{equation}
where$\ v=\kappa^{\alpha}\ $and$\ \Psi=v^{2}+v_{\theta}^{2}.$
\end{lem}
\begin{pf}
The proof is analogous to the proof of Lemma I1.12 in
Andrews \cite{And1} and we include it here for the convenience of a reader. Fix a $t>0$. Suppose that at $(\theta_0,t_0)\in
{I}\times[0,t]$ we have $\Psi(\theta_0,t_0) =
\sup_{{I}\times[0,t]}\Psi$. We may assume $t_0>0$ (otherwise we
are done). At the maximum of $\Psi$ we have $v_{\theta} (v + v_{\theta\theta}) = 0$. If the maximum of $\Psi$ is so large that $v_{\theta}^2 + v^2 > \sup_{I\times [0,t]}   v^2$, then we also have $v_{\theta}(\theta_0,t_0) \neq 0$.  This implies $v_{\theta\theta} + v = 0$ at $(\theta_0,t_0)$. Using that, a simple computation as in \cite{And1} shows that at the maximum point of $\Psi$ we have
\[
\begin{split}
\frac{\partial\Psi}{\partial t} &=  \alpha \kappa^{\alpha+1}\Psi_{\theta\theta}
-2\alpha^2(\alpha+1)\lambda(t)\kappa^{2\alpha-1}(\kappa_{\theta})^2-2\alpha\lambda(t)\kappa^{2\alpha+1} \\
&\le -2\alpha^2(\alpha+1)\lambda(t)\kappa^{2\alpha-1}(\kappa_{\theta})^2-2\alpha\lambda(t)\kappa^{2\alpha+1} \le 0,
\end{split}
\]
that is, the maximum of $\Psi$ is nonincreasing, which concludes the proof of the Lemma.
\end{pf}

\vskip 5 pt

Based on the gradient estimate, we obtain the following lemma.
\begin{lem} \label{blow_up_set}Assume that $\kappa(\theta,t)$ is the curvature
function of the evolving curves under the AP flow  or the LP flow,
and
$$\kappa(\theta_0,t_0)=\max_{I\times[0,t_0]}\kappa(\theta,t)$$
 for some $(\theta_0,t_0)\in I\times(0,T_{\max})$.
Then for any small $\varepsilon>0$, there exists a number $\delta>0$, depending only
on $\varepsilon$, such that
$$(1-\varepsilon)\kappa^\alpha(\theta_0,t_0) \leq \kappa^\alpha(\theta,t_0) + \epsilon\sqrt{C}$$
for all $\theta\in(\theta_0-\delta, \theta_0+\delta)$, where
$C$ is the constant only depending on the initial curve.
\end{lem}
\begin{pf} We adopt the notation $v=\kappa^\alpha$.
An easy integration gives that
\bess
v(\theta_0,t_0) &=& v(\theta,t_0) + \int_{\theta}^{\theta_0} v_{\theta}(\theta,t_0)\,d\theta\\
&\leq& v(\theta,t_0) + |\theta_0 - \theta| \max_{\theta\in I}|v_\theta(\theta,t_0)|.
\eess
Then by Lemma \ref{grad_est} we have
\bess
v(\theta_0,t_0) &\leq& v(\theta,t_0) + |\theta_0 - \theta| \sqrt{\max_{I\times[0,t_0]}v^2(\theta,t) + C}\\
&=& v(\theta,t_0) + |\theta - \theta_0| \sqrt{v^2(\theta_0,t_0) + C}\\
&\leq& v(\theta,t_0) + \delta
v(\theta_0,t_0)+ \delta \sqrt{C}. \eess Take $\delta:=\varepsilon$ and the lemma is proved.
\end{pf}

\vskip 10 pt

The lemma that follows will be useful to show the convergence of the flow, once the time-independent a priori estimate for the curvature is established.

\begin{lem}\label{priori_conv}
If there is a constant $C$ independent of time, such that
$$\max_{\theta\in I}\kappa(\theta,t)\leq C, \ \ t\in[0,T_{\max}),$$
with $T_{\max}$ being the maximal existence time, then the flow must exist for all time and converge smoothly to an $m$-fold circle as time goes to infinity.
\end{lem}
\begin{pf}
The proof is very similar to the one for the convergence of a nonlocal flow of embedded curves in Section 2.4 of \cite{TW1}, which uses the Lyapunov functional method to show desired convergence. The details are omitted here.
\end{pf}

\renewcommand{\theequation}{\thesection.\arabic{equation}}
\setcounter{equation}{0} \setcounter{prop}{0}
\renewcommand{\theprop}{\arabic{section}.\arabic{prop}}
\section{The evolution of area-preserving flow}

\vskip 10 pt

The main goal of this section is to prove Theorem \ref{them1}.
In a view of Lemma \ref{priori_conv}, in order to show the convergence of the AP flow  in certain cases and prove Theorem \ref{them1}  we need to show uniform curvature bounds along the flow.

\subsection{The convergence of the AP flow}


In this section, assuming global existence of the AP flow we show the integral $\int_I\kappa^\alpha\,d\theta$ is bounded for all times, and then by Lemma \ref{blow_up_set} we obtain the $L^\infty$-estimate of the curvature $\kappa$ along the time sequence $\{t_j\}_{j=1}^\infty$ such that $\max_{\theta\in I}\kappa(\theta,t_j)=\max_{I\times[0,t_j]}\kappa(\theta,t)$.

\begin{lem}\label{area_priori_estimate1}
If the AP flow exists for all times and if the $\lim_{t\rightarrow \infty}L(t)=L_\infty>0$, then for any fixed $d_0>0$, there is a sequence $t_j \in[(j-1)d_0,jd_0]$,
 such that
\bes
L(t_j)\int_{I} \kappa^\alpha(\theta,t_j)\,d\theta -2m\pi\int_I \kappa^{\alpha-1}(\theta,t_j)\,d\theta \rightarrow  0,\ \ \mathrm{as}\ \ t_j\rightarrow \infty,\label{area_limit_eq1}
\ees
and
\bes
\max_{j\geq j_0}\int_I \kappa^\alpha(\theta,t_j)\,d\theta \leq (2m\pi)^{\alpha+1}L_\infty^{-\alpha}, \label{bdint_area}
\ees
for some $j_0\geq 1$.
\end{lem}

\begin{pf} By the evolution equation of $L(t)$, we have
$$-\int_0^t\int_{X(\cdot,t)}\kappa(v-\lambda(\tau))\,dsd\tau = L(t)-L_0.$$
If the flow exists for all times, then
$$\int_0^{\infty}\int_{X(\cdot,t)}\kappa(v-\lambda(\tau))\,dsd\tau \leq L_0.$$
Denote
$$f(t) = \int_{X(\cdot,t)}\kappa(v-\lambda(t))\,ds.$$
The integral $\int_0^\infty f(\tau)\,d\tau$ is finite and thus we have
 $$\lim_{j\rightarrow\infty} \int_{(j-1)d_0}^\infty f(\tau)\,d\tau = 0.$$
 Noticing that $f(t)\geq 0$ for all times $t$, we have
 $$\lim_{j\rightarrow\infty} \int_{(j-1)d_0}^{jd_0} f(\tau)\,d\tau = 0.$$
By the mean value theorem, we can conclude that there exists a sequence $\{t_j\}_{j=1}^\infty$ with $t_j\in[(j-1)d_0,jd_0]$ such that
$$f(t_j) \rightarrow 0, \ \ j\rightarrow\infty,$$
that is,
$$\int_{I}\Big(\kappa^\alpha(\theta,t_j)-\frac{\int_I \kappa^{\alpha-1}(\theta,t_j)\,d\theta}{L(t_j)}\Big)\,d\theta \rightarrow  0,\ \ \mathrm{as}\ \ t_j\rightarrow \infty.$$
Hence, for any given $\epsilon>0$, there exists a $j_0$ depending on $\epsilon$, such that
\bes
\int_{I}\Big(\kappa^\alpha(\theta,t_j)-\frac{\int_I \kappa^{\alpha-1}(\theta,t_j)\,d\theta}{L(t_j)}\Big)\,d\theta \leq \epsilon,\ \  \forall\,j \geq j_0.\label{def_lim}
\ees

When $\alpha >1$, one may notice that
\bes
\int_I \kappa^{\alpha-1}\,d\theta \leq (2m\pi)^{1/\alpha}\Big(\int_I \kappa^\alpha\,d\theta\Big)^{(\alpha-1)/\alpha}.\label{area_glob_ineq1}
\ees
Taking $\epsilon = 2^{\alpha-1}(2m\pi)L_\infty^\alpha$ in (\ref{def_lim}), we can conclude that there is a $j_1$ such that
\bes
\int_{I}\Big(\kappa^\alpha(\theta,t_j)-\frac{\int_I \kappa^{\alpha-1}(\theta,t_j)\,d\theta}{L(t_j)}\Big)\,d\theta \leq 2^{\alpha-1}(2m\pi)L_\infty^\alpha,\ \  \forall\,j \geq j_1.\label{def_lim1}
\ees
Then we can claim from (\ref{def_lim1}) that
\bes
\int_I \kappa^\alpha(\theta,t_j)\,d\theta \leq 2^\alpha(2m\pi) L_\infty^{-\alpha},\ \  \forall\,j \geq j_1.
\label{glob_area_res2}\ees
Otherwise, if there exists a $j^*\geq j_1$ such that
$$\int_I \kappa^\alpha(\theta,t_{j^*})\,d\theta > 2^\alpha(2m\pi) L_\infty^{-\alpha},$$
then from (\ref{area_glob_ineq1}) we have
\bess
&&\int_{I}\Big(\kappa^\alpha(\theta,t_{j^*})-\frac{\int_I \kappa^{\alpha-1}(\theta,t_{j^*})\,d\theta}{L(t_{j^*})}\Big)\,d\theta\\
&\geq& \Big(\int_I\kappa^\alpha(\theta,t_{j^*})\,d\theta\Big)^{(\alpha-1)/\alpha}\Big[
\Big(\int_I\kappa^\alpha(\theta,t_{j^*})\,d\theta\Big)^{1/\alpha}
-\frac{(2m\pi)^{1/\alpha}}{L_\infty}\Big]\\
&>& 2^{\alpha-1}(2m\pi)L_\infty^\alpha,
\eess
which is a contradiction to (\ref{def_lim1}).

When $0<\alpha\leq 1$, noticing that
\[
\begin{split}
\int_I \kappa^{\alpha-1}\,d\theta  &= \int_I (\kappa^{-1})^{1-\alpha}\,d\theta
\le \left(\int_I \kappa^{-1}\, d\theta\right)^{1-\alpha} (2m\pi)^{\alpha} \\
&\leq  (2m\pi)^\alpha L^{1-\alpha},
\end{split}
\]
and $L(t) \geq L_\infty > 0$, we take $\epsilon = (2m\pi)^\alpha L_\infty^{-\alpha}$ in (\ref{def_lim}) to conclude that there is a $j_2$ such that
\bes
\int_I \kappa^\alpha(\theta,t_j)\,d\theta \leq 2(2m\pi)^\alpha L_\infty^{-\alpha},\ \  \forall\,j \geq j_2.\label{glob_area_res1}
\ees
Estimates (\ref{glob_area_res1}) for $0 < \alpha \le 1$  and
(\ref{glob_area_res2}) for $\alpha > 1$ conclude the proof of Lemma.
\end{pf}

\vskip 10 pt

We have showed the bound (\ref{bdint_area}) holds along a sequence. In order to prove the bound holds for all sufficiently large times we need to consider the evolution of $\int_I \kappa^\alpha\,d\theta$. Define
$$F(t) := \int_I \kappa^\alpha(\theta,t)\,d\theta,$$
or equivalently,
$$F(t) := \int_I \kappa^{\alpha+1}(s,t)\,ds.$$

\begin{lem}\label{area_evo_int}
For the function $F(t)$ defined as above, we have
\bes
F'(t)\leq C[F(t)^{\frac{4(2\alpha + 1)}{5\alpha + 1}}+F(t)^{2}],
\label{area_int_evo_ode}\ees
where the constant $C$ only depends on the exponents appearing in the Gagliardo-Nirenberg interpolation inequalities.
\end{lem}
\begin{pf}
A direct computation shows that
\bes
\frac{dF}{dt}
&=& \alpha \int k^{\alpha - 1} k_t d\theta = \alpha \int k^{\alpha} k_t ds\nonumber\\
&=& \alpha \int k^{\alpha} [(k^{\alpha})_{ss} + k^2(k^{\alpha} - \lambda(t))]\, ds\nonumber\\
&\le& -\alpha \int (k^{\alpha})_s^2 ds + \alpha \int k^{2\alpha + 2} ds. \label{evo_integ_area}
\ees
We shall use the Gagliardo-Nirenberg interpolation inequalities (called GN inequality for simplicity, see \cite{Ni1}): For a  function $u$ defined on $[0,L] $, which is sufficiently smooth, we have
$$||u^{(j)}||_{L^r}\leq C[||u||_{L^p}^{1-\theta}||u^{(k)}||_{L^q}^\theta + \|u\|_l], \ \ \theta\in(0,1),$$
where $r,p,q,j,k$ and $\theta$ satisfy $p,q,r>1, j\geq 0$,
$$\frac{1}{r} = j+\theta\Big(\frac{1}{q} - k\Big)+ (1-\theta)\frac{1}{p},$$
with
$$\frac{j}{k}\leq \theta \leq 1,$$
and $l>0$.
Here the constant $C$ depends on $r,p,q,j,k$ and $l$ only.

We bound the integral $\int_I \kappa^{2\alpha+2}\,ds$. Set
$$v=\kappa^\alpha.$$
By choosing $j=0,r=2+\frac{2}{\alpha},k=1,p=1+\frac{1}{\alpha},q=2$ and $l=1+\frac{1}{\alpha}$ in the GN inequality we have
\bess
\|v\|_{2+\frac{2}{\alpha}}\leq C_1[\|v\|_{1+\frac{1}{\alpha}}^{1-\theta}\|v_s\|_2^\theta + \|v\|_{1+\frac{1}{\alpha}}],
\eess
and thus
 \bess
\|v\|_{2+\frac{2}{\alpha}}^{2+\frac{2}{\alpha}}\leq C_2[\|v\|_{1+\frac{1}{\alpha}}^{(1-\theta)(2+\frac{2}{\alpha})}\|v_s\|_2^{\theta(2+\frac{2}{\alpha})} + \|v\|_{1+\frac{1}{\alpha}}^{2+\frac{2}{\alpha}}],
\eess
where
$$\theta=\frac{\alpha}{3\alpha+1}.$$
Then we use Young's inequality to obtain
 \bes
\|v\|_{2+\frac{2}{\alpha}}^{2+\frac{2}{\alpha}} &\leq& \frac{1}{2}\|v_s\|_2^2 + C_3\|v\|_{1+\frac{1}{\alpha}}^{\frac{4(2\alpha+1)(\alpha+1)}{\alpha (5\alpha+1)}} + C_2 \|v\|_{1+\frac{1}{\alpha}}^{2+\frac{2}{\alpha}}.\label{est_GN_area1}
\ees
Substituting (\ref{est_GN_area1})  into (\ref{evo_integ_area}), we have
\bess
\frac{d}{dt}\int_I\kappa^{\alpha+1}\,ds \leq C_3\|v\|_{1+\frac{1}{\alpha}}^{\frac{4(2\alpha+1)(\alpha+1)}{\alpha (5\alpha+1)}} + C_2 \|v\|_{1+\frac{1}{\alpha}}^{2+\frac{2}{\alpha}}.
\eess
\end{pf}

\vskip 10 pt
Now, from Lemma \ref{area_priori_estimate1} and Lemma \ref{area_evo_int}, we can obtain the bound on $\int_I \kappa^\alpha\,d\theta$ for all times. More precisely, we have the following Lemma.

\begin{lem}\label{glb_bd_int_area}
Under the assumptions of Lemma \ref{area_priori_estimate1},
there exists a time $T_0$, and a constant $C$ only depending on the limit of length $L_\infty$,  such that
\bes
\int_I\kappa^\alpha\,d\theta \leq C, \ \ \forall\,t\geq T_0.
\label{uni_bd_int}\ees
\end{lem}
\begin{pf}
Set $C_0=(2m\pi)^{\alpha+1}L_\infty^{-\alpha}$ to be a uniform constant that appears in the statement of Lemma \ref{area_priori_estimate1}. Since the integral $F(t)$ satisfies the ODE (\ref{area_int_evo_ode}), if the initial data $F(0)\leq C_0$, then there exists a $\delta_0>0$, such that
$$F(t)\leq 2C_0, \ \ \forall\,t\in[0,\delta_0].$$
Choose $d_0=\frac{\delta_0}{2}$, and then by Lemma \ref{area_priori_estimate1} find a sequence $\{t_j\}_{j=1}^\infty$ with $t_j\in[(j-1)d_0, jd_0)$ such that
$$\int_I\kappa^\alpha(\theta,t_j)\,d\theta \leq C_0, \ \ \forall\,j\geq j_0.$$
Using the above obeservation for solutions to ODE, we know that
$$\int_I\kappa^\alpha(\theta,t)\,d\theta \leq 2C_0, \ \ \mathrm{for}\ t\in[t_j,t_j+2d_0], \ \ \forall
\,j\geq j_0,$$
which implies that
$$\int_I\kappa^\alpha(\theta,t)\,d\theta \leq 2C_0, \ \ \forall\,t\geq t_{j_0}.$$
The proof is finished.
\end{pf}

\vskip 10 pt

Finally we can prove the convergence result for the AP flow. More precisely, we have the following Lemma.

\begin{lem}\label{convergence_AP}
Assume the initial curve $X_0$ is locally convex. If the AP flow exists for all
times and $\lim_{t\rightarrow\infty}L(t)=L(\infty)>0$, then $X(\cdot,t)$ converges smoothly to an $m$-fold circle.
\end{lem}
\begin{pf} In view of Lemma \ref{priori_conv}, we only need to show that the curvature of evolving curves has a time-independent upper bound. Indeed, if the claim does not hold, there exists a sequence $\{\theta_j\}_{j=1}^\infty\subset I$ and a sequence $\{t_j\}_{j=1}^\infty \rightarrow \infty$, such that
$$\kappa(\theta_j,t_j) = \max_{I\times[0,t_j]} \kappa(\theta,t)$$
and
$$\kappa(\theta_j,t_j)\rightarrow \infty, \ \mathrm{as} \ j\rightarrow\infty.$$
Then by Lemma \ref{blow_up_set}, we have
$$\int_I \kappa^\alpha(\theta,t_j)\,d\theta\rightarrow\infty, \ \mathrm{as} \ j\rightarrow\infty,$$
 a contradiction with the obtained bound for $\int_I\kappa^\alpha\,d\theta$ in Lemma \ref{glb_bd_int_area}. This implies there exists a uniform constant $C$ so that
\[\max_I \kappa(\cdot,t) \le C, \qquad \mbox{for all} \,\,\,\, t \in [0, \infty).\]
\end{pf}

\vskip 10pt

\subsection{The AP flow for highly symmetric curves}

When $\alpha=1$, one can mimick the proof of Gage \cite{Ga2} to show  the global existence of the AP flow when $X_0\in\mathcal{H}_{m,n}$, see \cite{WK1}. For $\alpha\neq 1$, the method of Gage does not apply and hence, it is not easy to obtain the similar gradient estimate as that one in  Corollary 3.5 of \cite{Ga2}. Here, we employ an isoperimetric bound established in Lemma 7.2 of \cite{And2} by Andrews to achieve our goal.

\begin{lem}(Andrews \cite{And2})\label{spt_harnack} For any curve in $\mathcal{H}_{m,n}$, its support function $h(\theta)$ satisfies
$$\sup_I h(\theta) \leq C\inf_I h(\theta),$$
for some constant $C$ only depending on $m$ and $n$.
\end{lem}

Consider the AP flow starting at an immersed, locally convex curve $X_0 \in \mathcal{H}_{m,n}$.
Immediately, we have two-sided bound for the evolving curve's support function.

\begin{lem}\label{h_bd_AP}For $X_0\in\mathcal{H}_{m,n}$, the support function of evolving curves under the AP flow satisfies
$$2r_0 \leq h(\theta,t) \leq 2R_0, \ (\theta,t)\in I\times[0,T_{\max}),$$
for some time-independent positive constants $r_0$ and $R_0$.
\end{lem}
\begin{pf}
Since $L(t)$ is nonincreasing and $L^2(t)\geq 4\pi |A(t)|$ (see an isoperimetric inequality of Rado in \cite{Oss1}),
we have
\begin{equation}
\label{eq-A0L0}
2\sqrt{\pi |A_0|} \leq L(t) \leq L_0.
\end{equation}
Notice that $|A(t)| = \frac{1}{2}\int_{X(\cdot,t)}h\,ds$ and
$$\inf_I h(\theta,t) L(t) \leq \int_{X(\cdot,t)}h\,ds \leq \sup_I h(\theta,t) L(t).$$
So by \eqref{eq-A0L0} we have $\inf_I h(\theta,t) \leq 2|A_0|/L(t) \leq |A_0|/\sqrt{\pi |A_0|}$, and by the monotonicity of $L(t)$ we have $\sup_I h(\theta,t) \geq 2|A_0|/L(t) \ge 2|A_0|/L_0.$ Then the two-sided bound for $h$ follows from Lemma \ref{spt_harnack}.
\end{pf}

\vskip 5 pt

We will use the two-side bound of $h$ in the proof of Thoerem \ref{them1} to
establish the upper bound on $\kappa$. The method is originally from
\cite{Tso1}.

\subsection{The AP flow for  Abresch-Langer  type curves}

The properties of Abresch-Langer type curves guarantee that evolving curves have `good' shape and thus
the estimates for curvature are feasible. In the following, $\kappa(\theta,t)$
 and $h(\theta,t)$ denote, as before, the curvature function and the support function of
 $X(.,t)$, respectively,  and they evolve from an Abresch-Langer type curve $X_0$ under the AP flow.

We first prove two lemmas in order to get some information about the shape of evolving
 $X(.,t)$.
\begin{lem}\label{symm}
Let $X(\cdot,t)$ be the solution to the AP flow starting at an Abresch-Langer type curve $X_0$. Then we have the following.
\begin{enumerate}
\item[(a)]
Both,  $\kappa(\theta,t)$ and $h(\theta,t)$  are symmetric with respect to
$\theta = 0$ and $\theta = m\pi/n$, for all times of the existence of the flow.
\item[(b)]
For all times $t$, both, $\kappa(\cdot,t)$ and $h(\cdot,t)$ attain their maximum at
$\theta = 0$; $h_\theta(\theta,t)$ and $\kappa_{\theta}(\theta,t)$ are negative on $(0,
m\pi/n)$, and consequently, $\kappa(\theta,t)$ and $h(\theta,t)$ are strictly decreasing in $(0,m\pi/n)$.
\end{enumerate}
\end{lem}
\begin{pf} It is easy to observe that (a) holds. We only show (b).
By differentiating the equation in (\ref{eqn2.2}), we see that the
function $w=h_\theta$ satisfies a parabolic equation \bess w_t =
a(\theta, t)w_{\theta\theta} + b(\theta,t)w, \ \ (\theta,t)\in
[-m\pi/n,m\pi/n]\times[0,T_{\max}) \eess where $ a(\theta,t) = b(\theta,t)= \alpha
\kappa^{\alpha+1}$. According to the Sturm comparison principle
(see \cite{Ang2} or \cite{Ma1}), the number of zeros of $w$ is
non-increasing in time. Since at initial time the function
$$w(\theta, 0) = \frac{\partial}{\partial\theta}h_{0}(\theta)$$ has exactly 2
zeros in $[-m\pi/n,m\pi/n]$ (a circle) by Property ($\mathcal{P}$),  the number of zeros of $w(\theta, t)$
cannot exceed two for all $t\in[0,T_{\max})$. On the other hand, the
reflectional symmetry of equation \eqref{eqn2.2} with respect to the
axis $\theta=0$ and $\theta=m\pi/n$ guarantees that $w(\theta,t)$
must vanish at $\theta=0$ and $m\pi/n$ for every $t\in[0,T)$. This
implies that $w(\theta,t)$ does not change its sign on $(-m\pi/n,0)$
and $(0,m\pi/n)$ for all $t\in[0,T_{\max})$. Then the conclusion for $h$
follows. The conclusion for $\kappa$ can be proved similarly.
\end{pf}

\vskip 5 pt

\begin{lem}\label{bdd_h_AL}
We have
$$h_0(m\pi/n)\leq h(\theta,t) \leq h_0(0), \ \ (\theta, t)\in I\times[0,T_{\max}).$$
\end{lem}
\begin{pf}
We claim that for any time $t\in(0,T_{\max})$
$$h_t<0 \ \ {\mathrm{at}}\ \theta=0;\  h_t>0 \ \ {\mathrm{at}}\ \theta=m\pi/n.$$
Indeed, since $\kappa(m\pi/n,t) \leq \kappa(\theta,t) \leq \kappa(0,t)$ by Lemma
\ref{symm}, we have
$${\kappa^\alpha(m\pi/n,t)} <{\dashint_I\ }{\kappa^\alpha(\theta,t)}\,d\theta < {\kappa^\alpha(0,t)}.$$
Then the claim is true in view of the equation $h_t={\dashint_I\
}\kappa^\alpha\,d\theta - \kappa^\alpha$. The proof is done.
\end{pf}

\vskip 10 pt

\subsection{Proof of Theorem \ref{them1}}

An isoperimetric inequality of Rado \cite{Oss1} tells
that for any closed, immersed curve,
$$L^2 \geq 4\pi\Sigma |m_j|A_j,$$
and thus
$$L^2 \geq 4\pi|\Sigma m_jA_j|=4\pi |A(t)|,$$
where $m_j$ and $A_j$ are the winding number and the enclosed area (which is nonnegaive in the usual sense) of the $j$-th component of
the curve, respectively. Since the flow preserves the (algebraic) area of the curve, we have
$$ L^2(t) \geq 4\pi|A_0|.$$
So, if $A_0\neq 0$ and $T_{\max}=+\infty$, then $\lim_{t\rightarrow\infty}L(t)>0$.

Assume now $A_0 < 0$ as in part (\ref{part-one}) of Theorem \ref{them1}. Our goal is to show the singularity occurs in finite time in this case. Assume on a contrary, that the flow exists forever. By Lemma \ref{convergence_AP},
we then obtain the convergence of the flow to an $m$-fold circle, which contradicts the fact that the flow preserves the negative enclosed area. Thus, if $A_0<0$, then a singularity must happen at some finite time.

Assume next the initial curve satisfies
$$0 < L_0^2 < 4m\pi A_0,$$
and exists for all times.
Since $dL(t)/dt\leq 0$ and $dA(t)/dt\equiv 0$, we have $L_0 \geq L(\infty):=\lim_{t\rightarrow\infty}L(t)$
and $A_0 = A(\infty):=\lim_{t\rightarrow\infty}A(t)$.
Thus,
\bes
L^2(\infty) < 4m\pi A(\infty).\label{isop_AP}
\ees
Then according to Lemma \ref{convergence_AP}
the flow converges smoothly to an $m$-fold circle as $t\rightarrow\infty$,
which implies that
\bess
L^2(\infty) = 4m\pi A(\infty).
\eess
This contradicts (\ref{isop_AP}). Thus, the singularity must happen at some finite time
during the evolution of the flow.

\vskip 10 pt
We now consider part (\ref{part-two}) of the Theorem, that is, the case when $A_0=0$. If the maximal existence time  $T_{\max}<\infty$, then it is well known the curvature must blow up as $t\rightarrow T_{\max}$. If $T_{\max}=\infty$, the curvature must blow up as $t\to \infty$, otherwise by Lemma \ref{priori_conv} we have the convergence to an $m$-fold circle enclosing a nonzero algebraic area, which contradicts the assumption $0 = A_0 = A(t)$, for all times $t\in [0,\infty)$.
Moreover, in the case $T_{\max}=\infty$,  the flow must go to a point as $t\to \infty$. Suppose  this is not true. It means that $\lim_{t\rightarrow\infty}L(t) > 0$. By Lemma \ref{convergence_AP}, the flow would then converge to an $m$-fold circle, which contradicts our assumption $A_0=0$.
 \vskip 10pt

For part (\ref{part-three}) of the Theorem, assume $X_0\in \mathcal{H}_{m,n}$ and $n > 2m$.
Fix a $t\in(0,T_{\max})$.
 Consider the quantity $\Phi =
\kappa^\alpha/(h-r_0)$ where $h(\theta,t)$ is the support function of evolving curves under the AP flow and $r_0$ is a constant from Lemma \ref{h_bd_AP} (by the same Lemma we have $h(\theta,t) \geq 2r_0$, which makes function $\Phi$ well defined).  Let the maximum of $\Phi$ over $I\times[0,t]$
be attained at $(\theta_0, t_0)$, $t_0>0$. At the point $(\theta_0,
t_0)$, we have \bess \frac{\partial \Phi}{\partial \theta} = 0, \
\frac{\partial \Phi}{\partial t} \geq 0, \ \mathrm{and}\
\frac{\partial^2 \Phi}{\partial \theta^2} \leq 0. \eess Since \bess
0\leq \frac{\partial \Phi}{\partial t} &=& \alpha
\kappa^{\alpha+1}\Phi_{\theta\theta} + \frac{2\alpha \kappa^{\alpha+1}h_\theta
\Phi_\theta}{h - r_0} + \frac{(\alpha+1)\kappa^{2\alpha}}{(h-r_0)^2}
- \frac{r\alpha \kappa^{2\alpha+1}}{(h-r_0)^2}\\
&&
- \lambda(t)\Big(\frac{\alpha \kappa^{\alpha+1}}{h_0-r}+\frac{\kappa^\alpha}{(h-r_0)^2}\Big)\\
&\leq& \frac{(\alpha+1)\kappa^{2\alpha}}{(h_0-r)^2} - \frac{r\alpha \kappa^{2\alpha+1}}{(h_0-r)^2}\\
&\leq & -\Phi^2[ r^{1+1/\alpha}\alpha\Phi^{1/\alpha} - (\alpha+1)]
\eess (where the inequality $h-r_0\geq r_0
> 0$ is used), we deduce that
$$\Phi(\theta_0,t_0)\leq r_0^{-(\alpha+1)}(1+\alpha^{-1})^\alpha.$$
If the maximum of $\Phi$ is attained at the initial time, we have
$$\Phi\leq \max_{I}\Phi(\theta,0).$$
Hence,
$$\Phi\leq \max\bigg\{r_0^{-(\alpha+1)}(1+\alpha^{-1})^\alpha, \max_{I}\Phi(\theta,0)\bigg\}:=M.$$
It follows that \bess
\kappa&\leq& M^{1/\alpha}(h - r_0)^{1/\alpha}\\
&\leq& M^{1/\alpha}(2R_0-r_0)^{1/\alpha}, \eess
where $R_0$ is the same constant as in Lemma \ref{h_bd_AP}.
At last,  the convergence of the flow
 is just an immediate result of Lemma \ref{priori_conv}.

\vskip 10 pt

In the case of part (\ref{part-four}) of the Theorem, due to the two-sided bound obtained in Lemma
\ref{bdd_h_AL}, the time-independent upper bound for curvature can be deduced immediately by the same proof of part (\ref{part-three}) of Theorem. The convergence then follows from Lemma \ref{priori_conv}. $\hfill\Box$

\vskip 10 pt

\renewcommand{\theequation}{\thesection.\arabic{equation}}
\setcounter{equation}{0} \setcounter{prop}{0}
\renewcommand{\theprop}{\arabic{section}.\arabic{prop}}
\section{The evolution of length-preserving flow}

For the LP flow, we can follow the steps in Sections 4.1 to show the convergence of global flow and then deduce the sufficient conditions for the occurrence of singularity at a finite time.

\subsection{The convergence of global LP flow}

\begin{lem}\label{priori_estimate_LP}
For the LP flow, if the flow exists for all time, then for any fixed $d_0>0$, there is a sequence $t_j \in[(j-1)d_0,jd_0]$ such that
\bes
L_0\int_{I} \kappa^\alpha(\theta,t_j)\,d\theta -2m\pi\int_I \kappa^{\alpha-1}(\theta,t_j)\,d\theta \rightarrow  0,\ \ \mathrm{as}\ \ t_j\rightarrow \infty,\label{limit_eq1}
\ees
and
$$\max_{j\geq 1}\int_I \kappa^\alpha(\theta,t_j)\,d\theta \leq C$$
for some constant $C$ independent of time, where $L_0$ is the length of an initial curve.
\end{lem}
\begin{pf}
By the evolution equation of $A(t)$, we have
$$-\int_0^t\int_{X(\cdot,t)}(v-\lambda(\tau))\,dsd\tau = A(t)-A_0.$$
Since $L^2(t)\geq 4\pi A(t)$ and $A(t)$ is nondecreasing in $t$, the limit $\lim_{t\rightarrow T_{\max}}A(t)$ is finite. If $T_{\max}=\infty$, then
$$\int_0^{\infty}\int_{X(\cdot,t)}(v-\lambda(\tau))\,dsd\tau > -\infty.$$
Since
$$\int_{X(\cdot,t)}(v-\lambda(t))\,ds \leq 0,$$
we can argue as in the proof of Lemma \ref{area_priori_estimate1} to conclude that for any fixed $d_0>0$ there exists a sequence $\{t_j\}_{j=1}^\infty$ with $t_j\in[(j-1)d_0,jd_0]$ such that
$$\int_{X(\cdot,t_j)}(v-\lambda(t))\,ds \rightarrow  0,\ \ \mathrm{as}\ \ t_j\rightarrow \infty,$$
that is,
$$\int_{I}\kappa^{\alpha-1}(\theta,t_j)\,d\theta - \frac{L_0}{2m\pi}\int_I \kappa^{\alpha}(\theta,t_j)\,d\theta \rightarrow  0,\ \ \mathrm{as}\ \ t_j\rightarrow \infty.$$
When $\alpha >1$, noticing that
$$\int_I \kappa^{\alpha-1}\,d\theta \leq (2m\pi)^{1/\alpha}\Big(\int_I \kappa^\alpha\,d\theta\Big)^{(\alpha-1)/\alpha},$$
we could employ the similar argument as in the proof of Lemma \ref{area_priori_estimate1} to show that $\max_{j\geq 1}\int_I \kappa^\alpha(\theta,t_j)\,d\theta \leq C_1$
for a constant $C_1$.
When $0<\alpha\leq 1$, noticing that
$$\int_I \kappa^{\alpha-1}\,d\theta  = \int_I (\kappa^{-1})^{1-\alpha}\,d\theta
\leq  (2m\pi)^\alpha L^{1-\alpha}$$
and $L(t) \equiv L_0$, it is easy to find a constant
$C_2$ such that
$$\max_{j\geq 1}\int_I \kappa^\alpha(\theta,t_j)\,d\theta \leq C_2.$$
This concludes the proof.
\end{pf}

\vskip 10 pt

Now, we can show the convergence of the LP flow if it
exists for all times.
\begin{lem}\label{convergence_LP}
Assume the initial curve $X_0$ is locally convex with the length $L_0>0$. If the LP flow exists for all
times, then $X(\cdot,t)$ converges smoothly to an $m$-fold circle.
\end{lem}
\begin{pf}  By checking the initial step in the proof of Lemma \ref{area_evo_int}, one can immediately observe that the ODE inequality (\ref{area_int_evo_ode}) also holds along the LP flow. Then arguing as in the proof of Lemma \ref{glb_bd_int_area}, we have the uniform bound for the integral $\int_I\kappa^\alpha\,d\theta$, for all times. Recall that the  estimate in Lemma \ref{blow_up_set} also holds for the LP flow. Following the proof of Lemma \ref{convergence_AP}, we use this integral estimate to get the uniform upper bound estimate for $\max_I\kappa(\cdot,t)$. The convergence then follows from Lemma \ref{priori_conv}.\end{pf}


\subsection{Proof of Theorem \ref{them2}}

Define $\overline{f} = \frac{1}{2m\pi}\int_I f\,d\theta$.
For the function $v(\theta,t)=\kappa^\alpha(\theta,t)$ we have
$$\overline{v} = \frac{1}{2m\pi}\int_I v\,d\theta = \lambda(t).$$
Recall that
$$E(t) = \int_I(v_{\theta})^2\,d\theta - \int_I(v-\overline{v})^2\,d\theta,$$
or
$$E(t) = \int_I (v_\theta)^2\,d\theta - \int_I v^2\,d\theta +
\frac{1}{2m\pi}\Big(\int_I v\,d\theta\Big)^2.$$

\begin{lem}\label{energy_LP}
For the energy $E(t)$ defined as above, we have
$$\frac{dE(t)}{dt} \leq 0.$$
\end{lem}
\begin{pf}
From the equation (\ref{curv_eqn_2}), we have
\bess
\int_I \frac{(v_t)^2}{\alpha v^p}\,d\theta &=& \int_I (v_{\theta\theta} + v -\overline{v})v_t\,d\theta \\
&=& -\frac{1}{2}\frac{d}{dt}\int_I [(v_\theta)^2 - v^2]\,d\theta - \overline{v}\int_I v_t\,d\theta,
\eess
where
\bess
\overline{v}\int_I v_t\,d\theta
= \frac{1}{4m\pi}\frac{d}{dt}\Big(\int_I v\,d\theta\Big)^2.
\eess
Thus,
\bess
-\frac{1}{2} \frac{dE(t)}{dt} = \int_I \frac{(v_t)^2}{\alpha v^p}\,d\theta \geq 0,
\eess
and the Lemma is proved.
\end{pf}

\vskip 10pt
One may ask what happens if the condition $E(0)<0$ does not hold for initial curve.
 A large class of rotationally symmetric curves belong to this case.
In fact, the Poincar$\acute{\mathrm{e}}$ inequality (or Wirtinger inequality) tells us the following is true.

\begin{lem} \label{main_cond2} If the initial curve is  locally convex, closed and has
  total curvature of $2m\pi$ and $n$-fold rotational symmetry with $m/n\leq 1$, then its curvature
  $k_0(\theta)$ satisfies
   \bes
   \int_I (v_0 - \overline{v_0})^2\,d\theta \leq \Big(\frac{m}{n}\Big)^2\int_I (v_{0\theta})^2\,d\theta.\label{poincare_ineq}\ees
\end{lem}

We are ready now to prove Theorem \ref{them2}.

\vskip 10 pt

\noindent{\textbf {Proof of Theorem \ref{them2}}.}
To prove part (1) of the Theorem, we argue by contradiction. Assume that for the initial curve we have $E(0)<0$, but that the flow exists for all time.  Then  Lemma \ref{convergence_LP} tells us that the flow must converge to an $m$-fold circle, whose energy is 0.
In view of the monotonicity of $E(t)$ in Lemma \ref{energy_LP}, we have $E(t)\geq 0$ for all $t\geq 0$, which contradicts the
assumption $E(0)<0$. Hence, the singularity must occur at some finite time.

If $E(0)=0$ and $k_0\not\equiv$ constant, we claim that
$(v_0)_{\theta\theta}+v_0-\overline{v_0}\neq 0$ must hold at some point of $I$. Indeed, if
$(v_0)_{\theta\theta}+v_0-\overline{v_0} = 0$ holds everywhere on $I$, it means that $v_0$ is a stationary solution of \eqref{curv_eqn_2}. Since $v_0\not\equiv$ constant, the initial curve is a non-circle locally convex closed curve, which either produces a singularity during the evolution, or exists globally and converges to an $m$-fold circle smoothly, according to Lemma \ref{convergence_LP}.  This contradicts the fact that $v_0$ is a noncircle stationary solution.  Hence, it holds that $(v_0)_{\theta\theta}+v_0-\overline{v_0}\neq 0$ at some point in $I$. Then by the continuity,
the same is true for some subinterval of $I$. By recalling the proof of Lemma \ref{energy_LP}, for $t$ close to 0, we have
\bess
 \frac{dE(t)}{dt} = -2\int_I \frac{(v_t)^2}{\alpha v^p}\,d\theta < 0,
\eess
which implies that $E(t)<0$ for $t>0$. Taking any $t_0 > 0$ to be the initial time, the above argument (since $E(t_0) < 0$) shows a singularity must happen at some finite time.

Next we show that if the initial curve satisfies $L_0^2 < 4m\pi A_0$, then a singularity must also occur at some finite time. Indeed, if we assume that the flow exists globally, then from Lemma \ref{convergence_LP} we have a smooth convergence to an $m$-fold circle. This implies that $L^2(\infty)=4m\pi A(\infty)$, a contradiction due to the monotonicity of $A(t)$.

\vskip 10 pt

To show part (2) of Theorem \ref{them2} we first consider the case when $\alpha>1$, or equivalently, $1<p<2$. By  equation (\ref{curv_eqn_2}) and integration by parts
we have
\bess
\frac{1}{\alpha(2-p)}\frac{d}{dt}\int_Iv^{2-p}\,d\theta &=& \int_Iv(v_{\theta\theta}+v-\bar{v})\,d\theta
= -E(t).
\eess
By Lemma \ref{main_cond2}, we have $E(t)\geq 0$ for $t\in[0,T_{\max})$,  since the evolving curves are rotationally symmetric. Hence, we have $\frac{d}{dt}\int_Iv^{2-p}\,d\theta \leq 0$. This implies there exists a constant
$C_1$, depending only on the initial curve, such that
$
\int_Iv^{2-p}\,d\theta \leq C_1 $
 for all $t\in[0,T_{\max})$. We claim there exists a constant
$C_2$ independent of time, such that
\bes
\max_{\theta\in I}\kappa(\theta,t) \leq C_2 \label{curvature_bd},
 \ees
 for all $t\in[0,T_{\max})$. To show the claim we argue by contradiction. Assume on the contrary that (\ref{curvature_bd})
 does not hold. Then Lemma \ref{blow_up_set} and the fact that $\int_I v^{2-p}\, d\theta \le C_1$ for all $t\in [0,T_{\max})$ yield contradiction.
 After the priori estimate (\ref{curvature_bd}) is established, we obtain the flow's global
existence and its smooth convergence to an $m$-fold circle
 as time goes to infinity  by Lemma \ref{priori_conv}.


Now we consider $0<\alpha\leq 1$. The case $m=n\,(=1)$ means that the curve is embedded and convex and has been studied in \cite{MZ1}  and \cite{TW1}. We only need to consider the case $m/n<1$. Since $p\geq2$ in this case, the above argument cannot be applied and hence we need different idea. Notice that $v$ is $2m\pi/n$-periodic. Then the Wirtinger inequality gives
$$E(t) \geq \Big(1-\Big(\frac{m}{n}\Big)^2\Big)\int_I (v_\theta(\theta,t))^2\,d\theta, \ \forall\ t\in[0,T_{\max}).$$
Since $E(t)\leq E(0)$, we have the gradient estimate
\bes
\int_I (v_\theta(\theta,t))^2\,d\theta \leq C, \ \forall\ t\in[0,T_{\max}),
\ees
with $C$ only depending on initial data. Fix any time $t\in[0,T_{\max})$. For any $\theta_1,\theta_2\in I$, it holds that
\be
\label{eq-v12}
\begin{split}
v(\theta_1,t) - v(\theta_2,t) &= \int_{\theta_2}^{\theta_1} v_{\theta}(\theta,t)\,d\theta\\
&\leq |\theta_1-\theta_2|^{\frac{1}{2}} \Big[\int_{\theta_1}^{\theta_2} (v_{\theta}(\theta,t))^2\,d\theta\Big]^{1/2}\\
&\leq C_1 \Big[\int_I (v_{\theta}(\theta,t))^2\,d\theta\Big]^{1/2}\\
&\leq C_2.
\end{split}
\ee
If we choose $\theta_2=\theta_2(t)$ such that
$$L_0 = \int_I\frac{d\theta}{\kappa(\theta,t)} = \frac{2m\pi}{\kappa(\theta_2(t),t)},$$
that is $\kappa(\theta_2(t),t) = 2m\pi/L_0$, and for $\theta_1$ to be any $\theta\in I$, then \eqref{eq-v12} yields the time-independent estimate
$$v(\theta,t) \leq  (2m\pi/L_0)^{\alpha} + C_2, \ (\theta,t)\in I\times[0,T_{\max}).$$
The desired convergence follows from Lemma \ref{priori_conv}.
 $\hfill\Box$

\subsection*{Acknowledgments.}
The first author thanks the NSF support in DMS-1056387. The second author is supported by NCTS and MoST of Taiwan with grant number 105-2115-M-007-007-MY3. The third author is supported by the Fundamental Research Funds for the Central Universities 2242015R30012, the NSF of China 11101078 and the Natural Science Foundation of Jiangsu Province BK20161412.






\bibliographystyle{elsarticle-num}



\end{document}